\documentclass {article}
\usepackage{authblk}
\usepackage{etex}
\usepackage[utf8]{inputenc}
\usepackage[english]{babel}
\usepackage{amsmath}
\usepackage{amsthm}
\usepackage{amssymb}
\usepackage{sectsty}
\usepackage{titlesec}
\usepackage{color}
\usepackage[color,matrix,arrow]{xy}
\usepackage{amsgen}
\usepackage{amstext}
\usepackage{amsbsy}
\usepackage{amsopn}
\usepackage{amsfonts}
\usepackage{eepic}
\usepackage{graphicx}
\usepackage{epsf}
\usepackage{pstricks}
\usepackage{enumerate}
\usepackage{eqnarray}
\usepackage{faktor}
\xyoption{all}

\usepackage{pgf}
\usepackage{tikz-cd}
\usetikzlibrary{automata, arrows.meta, positioning,decorations.pathmorphing}

\newcommand{\footrecall}[1]{%
} 
\usetikzlibrary{snakes,shapes,arrows,automata}

\titleformat*{\section}{\large\bfseries}
\titleformat*{\subsection}{\normalsize \bfseries}

\newcommand{\N}{\mathbb{N}}
\newcommand{\Z}{\mathbb{Z}}

\newcommand{\Q}{\mathbb{Q}}
\newcommand{\Fix}{\text{Fix}}
\newcommand{\rank}{\text{rank}}

\newcommand{\End}{\text{End}}
\newcommand{\Ker}{\text{Ker}}
\newcommand{\Evper}{\text{EvPer}}

\newcommand{\ord}{\varphi\text{-ord}}
\newcommand{\tord}{\tilde\varphi\text{-ord}}
\newcommand{\pord}{\varphi\text{-pord}}
\newcommand{\tail}{\varphi\text{-tail}}
\newcommand{\spe}{\varphi\text{-sp}}
\newcommand{\tspe}{\tilde\varphi\text{-sp}}
\newcommand{\Per}{\text{Per}}

\newcommand{\Orb}{\text{Orb}}
\newcommand{\Aut}{\text{Aut}}

\newcommand{\Img}{\text{Im}}

\usepackage{scalerel}
\newcommand{\mc}{\mathcal}

\theoremstyle{definition}
\newtheorem{theorem}{Theorem}[section]
\newtheorem{corollary}[theorem]{Corollary}

\newtheorem{proposition}[theorem]{Proposition}
\newtheorem{question}[theorem]{Question}

\newtheorem{lemma}[theorem]{Lemma}

\newtheorem{remark}[theorem]{Remark}
 
\begin{document}
 
 
\title{Quantifying Brinkmann's problem: relative $\varphi$-order and $\varphi$-spectrum}
\author{Andr\'e Carvalho}
\affil{Center for Mathematics and Applications (NOVA Math), NOVA SST

 2829–516 Caparica, Portugal

andrecruzcarvalho@gmail.com}
 
\maketitle
\begin{abstract}
We prove that the stable image of an endomorphism of a virtually free group is computable. For an endomorphism $\varphi$, an element $x\in G$ and a subset $K\subseteq G$, we say that the relative $\varphi$-order of $g$ in $K$, $\ord_K(g)$, is the smallest nonnegative integer $k$ such that $g\varphi^k\in K$. We prove that the set of orders, which we call $\varphi$-spectrum, is computable in two extreme cases: when $K$ is a finite subset and when $K$ is a recognizable subset. The finite case is proved for virtually free groups and the recognizable case for finitely presented groups. The case of finitely generated virtually abelian groups and some variations of the problem are also discussed.
\end{abstract}

\section{Introduction}
For a group $G$, the problem of deciding, on input two elements $x,y\in G$ and an endomorphism $\varphi\in \End(G)$, whether there exists $n\in \N$ such that $x\varphi^n=y$ or not, was introduced by Peter Brinkmann in \cite{[Bri10]}, who solved it for automorphisms of the free group. We refer to this problem as \emph{Brinkmann's Problem} ($BrP(G)$). The conjugacy variation of this problem turned out to be particularly important in proving the conjugacy problem for free-by-cyclic groups in \cite{[BMMV06]}. In \cite{[Log22]}, Logan solved several variations of this problem for general endomorphisms and used them to solve the conjugacy problem in ascending HNN-extensions of the free group, generalizing the work in \cite{[BMMV06]}. In fact, it is proved in \cite{[CD23]} that Logan's results imply decidability of $BrP(G)$ for endomorphisms (not necessarily injective) of the free group and this was extended to virtually free groups in \cite{[Car22c]}. Kannan and Lipton had already solved  in \cite{[KL86]}  the problem of deciding whether, given an $n\times n$ matrix $Q$ of rational numbers and two vectors of rational numbers $x,y\in \Q^n$, there is a natural number $i\in \N$ such that $xQ^i=y$, which is more general than  Brinkmann's Problem for free-abelian groups.

 The \emph{Generalized Brinkmann's Problem} on $G$, $GBrP(G)$, consists on deciding, given $x\in G$, $K\subseteq G$ and $\varphi\in \End(G)$,  whether there is some $n\in \N$ such that $x\varphi^n\in K$ or not. 
 Not much is known about $GBrP(G)$, even when $G$ is a finitely generated free group. However, $GBrP(G)$ was solved in \cite{[Car22c]} for finitely generated abelian groups, when $K$ is a coset of a finitely generated subgroup of $G$ and $\varphi$ is an automorphism. 
 When we consider the problem with restrictions on the class of sets or on the class of endomorphisms admitted as input, we write the restriction as the index.
 Clearly, the \emph{simple} version of the problem, $BRP(G)$, corresponds to the case where $K$ is a singleton and, for any group $G$ and any subclass of endomorphisms $\mc E$, decidability of $BrP_{\mc E}(G)$ implies decidability of  $GBrP_{(\mc E,  Fin(G))}(G)$, where $Fin(G)$ denotes the class of finite subsets of $G$. Indeed, if we want to decide whether, given $\varphi\in \mc E$, $x\in G$ and $K=\{y_1,\ldots, y_m\}\subseteq G$, there is some $n\in \N$ such that $x\varphi^n\in K$, we simply have to decide if there is some $n\in \N$ such that $x\varphi^n=y_i$, for $i\in [m]$.

Let $G$ be a finitely generated group, $K\subseteq G$ be a subset of $G$,  $g\in G$ be an element and $\varphi\in \End(G)$ be an endomorphism. 
Inspired by the terminology in \cite{[DVZ22]}, we say that the \emph{relative $\varphi$-order of $g$ in $K$} $\ord_K(g)$, is the smallest nonnegative integer $k$ such that $g\varphi^k\in K$. If there is no such $k$, we say that   $\ord_K(g)=\infty$.  The \emph{$\varphi$-spectrum} of a subset $\spe(K)$ is the set of relative $\varphi$-orders of elements in $K$, i.e., $\spe(K)=\{\ord_K(g)\mid g\in G\}$. A \emph{$\varphi$-preorder} is the set of elements of a given relative $\varphi$-order in $K$. For a nonnegative integer  $n\in \N$, we denote the $\varphi$-preorder of $n$ in $K$ by $\pord_K(n)=\{g\in G\mid \ord_K(g)=n\}.$

 Naturally, if we can decide membership in $K$ and $GBrP(G)$, we can compute $\ord_K(g)$: if $GBrP(G)$, with input $g$, $K$ and $\varphi$ answers \texttt{NO}, then  $\ord_K(g)=\infty$; otherwise, we simply start computing the $\varphi$-orbit of $g$ and checking membership of $g\varphi^k$ in $K$ until we get a positive answer.
 
The goal of this paper is to, in some sense, quantify Brinkmann's problem, by showing that, in some cases, relative $\varphi$-orders and the $\varphi$-spectrum can be computed. Our results concern mainly the class of virtually free groups, but we also consider the problem in virtually abelian groups. As highlighted above, the generalized version of Brinkmann's problem can be quite challenging, being only known to be decidable in the class of finitely generated abelian groups. Since computability of orders implies in particular decidability of $GBrCP$, we will restrict ourselves to the study of two extreme cases: the case where $K$ is finite and the case where $K$ is recognizable, i.e., a union of cosets of a finite index subgroup. This, in particular, quantifies the simple version of the problem, $BrP(G).$ Some discussion on more general sets will be taken in the final section of the paper.

For an endomorphism $\varphi\in\End(G)$, the \emph{stable image of $\varphi$} is $$\varphi^\infty(G)=\bigcap_{i=1}^\infty \varphi^i(G).$$
This notion was introduced in \cite{[IT89]} and it has been a most useful tool in reducing problems concerning endomorphisms to ones concerning only automorphisms, since it was proved in \cite{[IT89]} that, if $G$ is a finitely generated free group, then $\varphi^\infty(G)$ is finitely generated and $\varphi|_{\varphi^\infty(G)}$ is an automorphism. Mutanguha proved in \cite{[Mut21]} that a basis for the stable image of an endomorphism of a finitely generated free group can be computed, which, combined with previous work by Bogopolski and Maslakova,  implies that $\Fix(\varphi)$ can be computed for general endomorphisms of a finitely generated free group. 

This concept appears naturally when we study the $\varphi$-spectrum of an endomorphism, so our first result is an extension of Mutanguha's to the class of finitely generated virtually free groups.

\newtheorem*{stable image computable}{Theorem \ref{stable image computable}}
\begin{stable image computable}
There exists an algorithm with input a finitely generated virtually free group $G$ and an endomorphism $\varphi$ of $G$ and output a finite set of generators for  $\varphi^\infty(G)$. 
\end{stable image computable}

With this at hand, we begin the study of the computation of  the $\varphi$-spectrum for endomorphisms of finitely generated virtually free groups and prove the following that it is computable for finite subsets of $G$.
\newtheorem*{main finite}{Theorem \ref{main finite}}
\begin{main finite}
There exists an algorithm with input a finitely generated virtually free group $G$, an endomorphism $\varphi$ of $G$ and a finite set $K=\{g_1,\ldots, g_k\}\subseteq G$ and output $\spe(K)$.
\end{main finite}

 We then prove that if one is able to prove an analogue of Theorem \ref{stable image computable} for finitely generated virtually abelian groups, an analogue of Theorem \ref{main finite} would follow. We then move from the finite case to the case where the subset $K$ is \emph{recognizable} and we can prove computability of orders, preorders and spectrum in the very general setting of finitely presented groups.
 \newtheorem*{recognizable computability}{Theorem \ref{recognizable computability}}
 \begin{recognizable computability}
Let $G$ be a finitely presented group, $K\in Rec(G)$, $g\in G$ and $\varphi\in \End(G)$. Then $\ord_K(g)$, $\spe(K)$ and the $\varphi$-preorders are computable.
\end{recognizable computability}

We then present some open questions and comments on variations of these problems.

The paper is organized as follows. In Section 2, we present some preliminaries on subsets of groups and on the $\varphi$-spectrum of an endomorphism.  In Section 3, we study the computability of the stable image of an endomorphism, proving it is computable in case the group is virtually free and we also discuss this problem in the abelian case. In Section 4, we study the computation of the $\varphi$-spectrum of an endomorphism for a finite set $K$, proving it is computable in the virtually free case and showing that the virtually abelian case depends on the computation of stable images. In Section 5, we consider the case where $K$ is a recognizable set and prove that, in this case, orders, preorders and the spectrum are all computable for endomorphisms of finitely presented groups. We end with some open questions and natural variations on this problem in Section 6. 

\section{Preliminaries}
Given a finitely generated group $G=\langle A\rangle$, a finite generating set $A$ and a set of formal inverses $A^{-1}$, write $\tilde A=A\cup A^{-1}$. There is a canonical (surjective) homomorphism $\pi:\tilde A^*\to G$
mapping $a\in \tilde A$ (resp. $a^{-1}\in \tilde A$) to  $a\in  G$ (resp. $a^{-1}\in  G$).

A subset $K\subseteq G$ is said to be \emph{rational} if there is some rational language $L\subseteq \tilde A^*$ such that $L\pi=K$.
The set of all such subsets is denoted by $\text{Rat } G$. Rational subsets generalize the notion of finitely generated subgroups.

\begin{theorem}\cite[Theorem III.2.7]{[Ber79]}
\label{anisimovseifert}
Let $H$ be a subgroup of a group $G$. Then $H\in \text{Rat } G$ if and only if $H$ is finitely generated.
\end{theorem}

A subset $K\subseteq G$ is said to be \emph{recognizable} if $K\pi^{-1}$ is rational.  The class of recognizable subsets of a finitely generated group $G$, $Rec(G)$, is effectively closed under Boolean operations and inverse morphisms. Recognizable subsets generalize the notion of finite index subgroups.

\begin{proposition}
\label{rec fi}
Let $H$ be a subgroup of a group $G$. Then $H\in  Rec(G)$ if and only if $H$ has finite index in $G$.
\end{proposition}

In fact, if $G$ is a group and $K$ is a subset of $G$ then $K$ is recognizable if and only if $K$ is a (finite) union of cosets of a subgroup of finite index.

\subsection{Virtually free and virtually abelian groups}

A subgroup $H$ of a group $G$ is \emph{fully invariant} if $\varphi(H)\subseteq H$ for every endomorphism $\varphi$ of $G$.

\begin{lemma}\cite[Lemma 4.1]{[Car22]}
\label{fullyinvariant}
Let $G$ be a group, $n$ be a natural number and $N$ be the intersection of all normal subgroups of $G$ of index $\leq n$. Then $N$ is fully invariant, and if $G$ is finitely generated, then $N$ has finite index in $G$.
\end{lemma}

  When we take a finitely generated virtually free group as input, we assume that we are given a decomposition as a disjoint union
\begin{align}
\label{decomp}
G=Fb_1\cup Fb_2\cup \cdots \cup Fb_m,
\end{align}
where $F=F_A\trianglelefteq G$ is a finitely generated free group, and a presentation of the form $\langle A,b_1,\ldots, b_m \mid R\rangle$, where the relations in $R$ are of the  form $b_ia=u_{ia}b_{i}$ and $b_{i}b_{j}=v_{ij}b_{r_{ij}}$, with  $u_{ia}, v_{ij} \in F_A$ and $r_{ij}\in [m]$, $i,j=1\ldots,m$, $a\in A$.
In the proof of \cite[Theorem 4.2]{[Car22]}, using Lemma \ref{fullyinvariant}, it is shown that given a finitely generated virtually free group, we can assume that the subgroup $F$ is a fully invariant normal subgroup of $G$.

Analogously, if we are given a f.g. virtually abelian group $G$, we can also assume that we are given a fully invariant free-abelian subgroup, together with a decomposition of the form (\ref{decomp}), since $G$ has a free-abelian subgroup of finite index and subgroups of free-abelian groups are again free-abelian.

Let $G$ be a finitely generated group, $K\subseteq G$ be a subset of $G$,  $g\in G$ be an element and $\varphi\in \End(G)$ be an endomorphism. 
Inspired by the terminology in \cite{[DVZ22]}, we say that the \emph{relative $\varphi$-order of $g$ in $K$} $\ord_K(g)$, is the smallest nonnegative integer $k$ such that $g\varphi^k\in K$. If there is no such $k$, we say that   $\ord_K(g)=\infty$.  The \emph{$\varphi$-spectrum} of a subset $\spe(K)$ is the set of relative $\varphi$-orders of elements in $K$, i.e., $\spe(K)=\{\ord_K(g)\mid g\in G\}$. A \emph{$\varphi$-preorder} is the set of elements of a given relative $\varphi$-order in $K$. For a nonnegative integer  $n\in \N$, we denote the $\varphi$-preorder of $n$ in $K$ by $\pord_K(n)=\{g\in G\mid \ord_K(g)=n\}.$

It is 
clear from the definitions that  
 $$\pord_K(n)=K\varphi^{-n}\setminus \bigcup_{i=0}^{n-1} K\varphi^{-i}
 $$
and 
 \begin{align}
 \label{defspecpord}
 \spe(K)=\{n\in \N\mid \pord_K(n)\neq \emptyset\}.
 \end{align}
 
  The following lemma is obvious by the definition.
  \begin{lemma}
 Let $G$ be a group, $K\subseteq G$, $\varphi\in \End(G)$ and $g\in G\setminus K$. 
Then, $$\ord_K(g\varphi)=\ord_K(g)-1.$$ 
 \end{lemma}
 
 For $n\in \N$, we denote the set $\{0,1,\cdots, n\}$ by $[n]_0$.
 \begin{corollary}
 Let $G$ be a group, $K\subseteq G$ and $\varphi\in \End(G)$. If, for some $n\in \N$, $\pord_K(n)=\emptyset$, then $\pord_K(m)=\emptyset$ for all $m\geq n$. 
 Moreover,  $$\spe(K)=
 \begin{cases}
 [n]_0\quad &\text{ if  $n=\min\{m\in \N\mid \pord_K(m+1)=\emptyset\}$}\\
 \N\quad &\text{ if  $\forall m\in \N, \, \pord_K(m)\neq\emptyset$}
 \end{cases}
 $$ 
 \end{corollary}
 \noindent\textit{Proof.} Suppose that $\pord_K(n)=\emptyset$. If  $\pord_K(n+1)\neq \emptyset$, let $g\in \pord_K(n+1)$. Then $\ord_K(g\varphi)=n$ and so $g\in \pord_K(n)$, which is absurd. So, $\pord_K(n)=\emptyset$ implies that $\pord_K(n+1)=\emptyset$ and the result follows by induction.

The observation about $\spe(K)$ can be proved in the same way. 
 \qed\\

  \section{Stable images}
  
 For $\varphi\in\End(G)$, the \emph{stable image of $\varphi$} is $$\varphi^\infty(G)=\bigcap_{i=1}^\infty \varphi^i(G).$$
This notion was introduced in \cite{[IT89]}, where it was proved that, if $G$ is a finitely generated free group, then $\varphi^\infty(G)$ is finitely generated and $\varphi|_{\varphi^\infty(G)}$ is an automorphism. 

Mutanguha proved in \cite{[Mut21]} that a basis for the stable image of an endomorphism of a finitely generated free group can be computed, which, combined with previous work by Bogopolski and Maslakova,  implies that $\Fix(\varphi)$ can be computed for general endomorphisms of a finitely generated free group.

Analogously, a similar result holds for free-abelian groups: if $G$ is a f.g. free-abelian group, then every subgroup of $G$ is finitely generated and so $\varphi^\infty(G)$ is finitely generated. Also, the restrictions $\varphi_k:\Img(\varphi^k)\to \Img(\varphi^{k+1})$ of $\varphi$ are necessarily surjective and, for large enough $k$, we have that it must be injective too, since it must eventually happen that $\rank(\Img(\varphi^k))=\rank(\Img(\varphi^{k+1}))$, and free-abelian groups are Hopfian.

Given a word $w\in X^*$ and a letter $c\in X$, we denote by $n_c(w)$ the number of occurrences of the letter $c$ in the word $w$.

 \begin{theorem}\label{stable image computable}
There exists an algorithm with input a finitely generated virtually free group $G$ and an endomorphism $\varphi$ of $G$ and output a finite set of generators for  $\varphi^\infty(G)$. 
 \end{theorem}
 \noindent\textit{Proof.} Write $G$ as a disjoint union of the form  $$G=Fb_1\cup \cdots \cup Fb_m,$$
 where $F$ is a fully invariant free normal subgroup of $G$.
 Then   
 \begin{align*}
 \varphi^\infty(G)&=\varphi^\infty(G)\cap G\\
& =(\varphi^\infty(G)\cap Fb_1)\cup \cdots \cup (\varphi^\infty(G)\cap Fb_m).
 \end{align*}
 We will prove that each of the subsets $\varphi^\infty(G)\cap Fb_i$ is rational and computable, which implies that  $\varphi^\infty(G)$ is finitely generated and computable, since it is a finite union of (computable) rational subsets (see Theorem \ref{anisimovseifert}).
 
 Fix $i\in [m]$. Since, for all $k\in \N$, $G\varphi^k=\bigcup_{j=1}^m F\varphi^k(b_j\varphi^k) $, we have that 
 \begin{align*}
 \varphi^\infty(G)\cap Fb_i &=\left(\bigcap_{k\in \N} G\varphi^k \right) \cap  Fb_i\\
 &=\bigcap_{k\in \N} (G\varphi^k   \cap Fb_i)\\
 &=\bigcap_{k\in \N} \bigcup_{j=1}^m (F\varphi^k(b_j\varphi^k)   \cap Fb_i).
 \end{align*}
 
 Since $F$ is fully invariant, we have that, for all $k\in \N$,  $(Fb_i)\varphi^k\subseteq F(b_i\varphi^k)$. Hence, the mapping $\theta: G/F\to G/F$ defined by $Fb_i\mapsto F(b_i\varphi)$ is a well-defined endomorphism.
 Since $|G/F|=m$ is finite, we can compute the orbit $\Orb_\theta(Fb_i)$ of $Fb_i$ through $\theta$. In particular, we can check if $Fb_i$ is periodic. 
 
 If it is not, then   $Fb_i\not\in \Img(\theta^k)$, for $k> m$.  This means that for $k> m$, and $j\in[m]$,
$$ F\varphi^k(b_j\varphi^k)   \cap Fb_i\subseteq F(b_j\varphi^k)   \cap Fb_i=(Fb_j)\theta^k\cap Fb_i=\emptyset,$$
since $Fb_i\not\in \Img(\theta^k)$ and so $(Fb_j)\theta^k$ and $Fb_i$ are distinct, thus disjoint, cosets. Hence, 
 $ \varphi^\infty(G)\cap Fb_i =\emptyset$.

 If $Fb_i$ is periodic, say of period $P$, then we claim that 
 \begin{align}
 \label{eq1}
 \varphi^\infty(G)\cap Fb_i=\bigcap_{k\in \N}(G\varphi^k\cap Fb_i)=\bigcap_{k>m}(G\varphi^k\cap Fb_i)=\bigcap_{k> m} F\varphi^{Pk}b_i\varphi^{Pk}.
 \end{align}
 
The first two equalities are obvious and it is clear that $$\bigcap_{k> m} F\varphi^{Pk}b_i\varphi^{Pk}\subseteq \bigcap_{k>m}(G\varphi^k\cap Fb_i)$$
since $Fb_i$ is a periodic point of $\theta$ with period $P$. Now we prove the reverse inclusion.
Let $x\in \bigcap_{k>m}(G\varphi^k\cap Fb_i)$ and fix $k>m$.
Then $x\in G\varphi^{2Pk}\cap Fb_i$, which means that there are some $f\in F$ and $j\in [m]$ such that $$(fb_j)\varphi^{2Pk}=(fb_j)\varphi^{Pk}\varphi^{Pk}=x\in Fb_i.$$
Since $k>m$, it follows that $(Fb_j)\theta^{Pk}$ must belong to the periodic orbit of $Fb_i$. Since it gets mapped to $Fb_i$ after $Pk$ applications of $\theta$, we must have that $(Fb_j)\theta^{Pk}=Fb_i$, and so $(fb_j)\varphi^{Pk}\in Fb_i$ and 
$$x=(fb_j)\varphi^{Pk}\varphi^{Pk}\in (Fb_i)\varphi^{Pk}=F\varphi^{Pk}b_i\varphi^{Pk}.$$
Since $k$ is arbitrary, we have proved (\ref{eq1}).

Now, since $Fb_i$ is a periodic point of $\theta$, this means that there is some $y\in F$ such that $b_i\varphi^P=yb_i$. Let $c$ be a letter not belonging to $F$, $F'=F*\langle c|\rangle$ and $\psi:F'\to F'$ be the endomorphism that applies $\varphi^P$ to the elements in $F$ and maps $c$ to $yc$. We will prove that 
\begin{align*}
\varphi^\infty(G)\cap Fb_i=\bigcap_{k> m} F\varphi^{Pk}b_i\varphi^{Pk}=\left(\bigcap_{k> m} F\psi^{k}c\psi^{k}\right)c^{-1}b_i=(\psi^\infty(F')\cap Fc)c^{-1}b_i.
\end{align*}
The first equality is simply (\ref{eq1}). The second equality follows from the fact that $$b_i\varphi^{Pk}=\left(\prod_{j=0}^{k-1}y\varphi^{(k-j-1)P}\right) b_i=(c\psi^{k})c^{-1}b_i.$$
Finally, it is clear that   $(\bigcap_{k> m} F\psi^{k}c\psi^{k})=\bigcap_{k> m} (Fc)\psi^{k}\subseteq (\psi^\infty(F')\cap Fc)$, and so $$(\bigcap_{k> m} F\psi^{k}c\psi^{k})c^{-1}b_i\subseteq (\psi^\infty(F')\cap Fc)c^{-1}b_i.$$
We only have to prove that $(\psi^\infty(F')\cap Fc)\subseteq \bigcap_{k> m} F\psi^{k}c\psi^{k} $. We start by proving that
\begin{align*}
\forall k\in \N,\; \forall x\in F'\; (x\psi^k\in Fc\implies x\psi^k\in (Fc)\psi^k)
\end{align*}
Let $k\in \N$, $x\in F'$ and $X$ be a basis for $F$. We proceed by induction on $|x|$. If $|x|=1$, then $x\psi^k\in Fc$ implies that $x=c$, since $F\psi^k\subseteq F$ and so $x\psi^k=c\psi^k\in (Fc)\psi^k$. So assume the claim holds for all $x$ such that $|x|\leq n$ and take $x$ such that $|x|=n+1$. 
Write $x=x_1\cdots x_{n+1}$ and $x\psi^k=w_1\cdots w_{n+1}$, where $x_i\in \widetilde{X\cup\{c\}}$ and $w_i=x_i\psi^k$. Let $w=w_1\cdots w_{n+1}\in (\widetilde{X\cup\{c\}})^*$. Then 
$n_c(x)=n_c(w)$ and $n_{c^{-1}}(w)=n_{c^{-1}}(x)$.  Choose any cancellation order on 
$w_1\cdots w_{n+1}$. We have that the result of cancellation, $\overline{w_1\cdots w_{n+1}}$, belongs to $Fc$. The $c$ that survives the cancellation must belong to some $w_i$ and it appears as the last letter of $w_i$, by definition of $\psi$. If $i\leq n$, then we have that $x\psi^k=\overline{w_1\cdots w_i}=(x_1\cdots x_i)\psi^k$, and so, by the induction hypothesis, it follows that $x\psi^k=(x_1\cdots x_i)\psi^k\in (Fc)\psi^k$.
If, on the other hand the surviving $c$ belongs to $w_{n+1}$, then either $n_c(w)=1$ and $n_{c^{-1}}(w)=0$, in which case the claim is obvious since $x\in Fc$, or $n_c(w)>1$ and $n_{c^{-1}}(w)=n_c(w)-1>0$ and there is cancellation between $c$'s and $c^{-1}$'s. Consider the first cancellation occurring between $c$'s and $c^{-1}$'s and let $i, j$ be such that the cancelled $c$ and the cancelled $c^{-1}$ belong to $w_i$ and $w_j$, respectively, and so $x_i=c$ and $x_j=c^{-1}$ and $c$ (resp. $c^{-1}$) is the last (resp. first) letter of $w_i$ (resp. $w_j$).  Suppose that $i<j$. Then $j>i+1$, since otherwise $x$ is not reduced. Thus, we have that $\overline{w_{i+1}\cdots w_{j-1}}=1$, i.e.,  ${x_{i+1}\cdots x_{j-1}}\in \Ker(\psi^k)$ and $x\psi^k=x_1\cdots x_ix_j\cdots x_{n+1}$.  By the induction hypothesis, it follows that $x\psi^k=(x_1\cdots x_ix_j\cdots x_{n+1})\psi^k\in (Fc)\psi^k$. If $i>j$, we proceed analogously. Again we have that $i>j+1$, since otherwise $x$ is not reduced. Since the first letter of $w_j$ cancels with the last letter of $w_i$, we have that   $\overline{w_{j}\cdots w_{i}}=1$, i.e.,  ${x_{j}\cdots x_{i}}\in \Ker(\psi^k)$ and $x\psi^k=x_1\cdots x_{j-1}x_{i+1}\cdots x_{n+1}$.  By the induction hypothesis, it follows that $x\psi^k=(x_1\cdots x_{j-1}x_{i+1}\cdots x_{n+1})\psi^k\in (Fc)\psi^k$.

Hence,  $\Img(\psi^k)\cap Fc \subseteq  (Fc)\psi^k$ and $$(\psi^\infty(F')\cap Fc)\subseteq \bigcap_{k\geq m} (Fc)\psi^{k}=\bigcap_{k\geq m} F\psi^{k}c\psi^{k}.$$

By \cite{[Mut21]}, $\psi^\infty(F')$ is finitely generated (thus rational) and computable, so $\varphi^\infty(G)\cap Fb_i=(\psi^\infty(F')\cap Fc)c^{-1}b_i$ is also rational and can be effectively computed.
 \qed\\

As far as the author knows, the same problem is not known to be decidable for finitely generated virtually abelian groups.
\begin{question}\label{question stable image v abelian}
Can we compute the stable image of an endomorphism of a finitely generated virtually abelian group?
\end{question}

It is likely that answering the question for endomorphisms of free-abelian groups would suffice to answer Question \ref{question stable image v abelian}, using a similar method to the one used to extend the result in the free case to finite extensions of free groups. However, that does not seem to be known.
\begin{question}\label{question stable image free abelian}
Can we compute the stable image of an endomorphism of a finitely generated free-abelian group? Equivalently,  letting $G=\Z^m$ and $Q\in \mc M_{m\times m}(\Z)$, can we compute $Q^\infty(G)$?
\end{question}
It is an easy observation that, to study the question above, we may assume injectivity of $Q$. 
 Consider the chain  $$\Img(Q)\supseteq \Img(Q^2) \supseteq \Img(Q^3)\supseteq \cdots $$
 and  the surjective homomorphisms $Q_i:\Img(Q^i)\to\Img(Q^{i+1})$ obtained by restriction of $Q$.
By Hopfianity, we have that, for $i\in \N$, $Q_i$ is injective if and only if $\rank(\Img(Q^i))=\rank(\Img(Q^{i+1}))$, which must happen for some computable $i\in \N$, since $\rank(\Img(Q^{i+1}))\leq \rank(\Img(Q^i))$.
 Fix such $i$. Then we consider the injective endomorphism $Q':\Img(Q^i)\to\Img(Q^i)$ given by $Q_i$. Clearly $Q'^\infty(\Img(Q^i))=Q^\infty(G)$, and so we may assume that our endomorphism $Q$ is injective.

\begin{remark}
\label{rem: phik injetivo}
In the proof of \cite[Proposition 5.5]{[Car22]}, the author shows that, for a f.g. virtually free group $G$ and an endomorphism $\varphi\in \End(G)$,  we can compute a constant $k$ such that the restriction of $\varphi$, $\varphi_k:\Img(\varphi^k)\to \Img(\varphi^{k+1})$ is injective. The same reasoning can be applied when $G$ is a f.g. virtually abelian group, as done above.
\end{remark}

We will denote the restriction of $\varphi$ to $\Img(\varphi^k)$ by $\varphi_k$ throughout the paper. Sometimes, we will also restrict the codomain to $\Img(\varphi^{k+1})$, or even $\Img(\varphi^k)$, but it should be clear from the situation what the codomain is.

\section{Finite $K$}
We will start by the simpler case where $K$ is a singleton. Computing relative $\varphi$-orders, which corresponds to solving $BrP(G)$, is known to be possible in case the group is f.g. virtually free or f.g.  abelian.

We now present two technical lemmas related to kernels of powers of $\varphi$ that will be useful later.

\begin{lemma} \label{lema nucleos1}
Let $G$ be a group, $\varphi\in \End(G)$, $i, j\in \N$ be such that $i\geq j$  and $a_1,\ldots, a_n\in \Ker(\varphi^i)$. The following are equivalent:
\begin{enumerate}
\item  $\Ker(\varphi^i)=a_1\Ker(\varphi^{j})\cup \cdots \cup a_n\Ker(\varphi^{j})$; 
\item $\Ker(\varphi_{j}^{i-j})= \{a_1\varphi^{j},\ldots, a_n\varphi^{j}\}$.
\end{enumerate}
\end{lemma}
\noindent\textit{Proof.}  Assume 1. Clearly, for all $r\in \{1,\ldots, n\}$, $a_r\varphi^{j}\varphi_{j}^{i-j}=a_r\varphi^i=1$.
 Now, let $x\in \Ker(\varphi_{j}^{i-j})$. Then, there is some $y\in G$ such that $x=y\varphi^{j}$ and $y\varphi^i=y\varphi^{j+i-j}=1$. Hence $y\in a_r\Ker(\varphi^{j})$ for some $r\in \{1,\ldots,n\}$, and so $x=y\varphi^{j}=a_r\varphi^{j}$.
 
Now we prove that  2 $\implies$ 1.  It is clear that   $a_1\Ker(\varphi^{j})\cup \cdots \cup a_n\Ker(\varphi^{j})\subseteq \Ker(\varphi^i)$. Now, let $x\in \Ker(\varphi^i)$. Then, 
 $1=x\varphi^i=x\varphi^{j}\varphi_{j}^{i-j}$ and so $x\varphi^{j}=a_r\varphi^{j}$ for some $r\in \{1,\ldots, n\}$  and so $x\in a_r\Ker(\varphi^{j})$.
  \qed\\

As a particular case, we obtain the following.

\begin{lemma} 
\label{lem: nucleosphik}
Let $G$ be a group, $k\in \N$ and $\varphi\in \End(G)$. The following are equivalent:
\begin{enumerate}
\item $\varphi_k$ is injective;
\item $\Ker(\varphi^k)= \Ker(\varphi^i)$ for all $i> k$;
\item $\Ker(\varphi^k)= \Ker(\varphi^i)$ for some $i> k$.
\end{enumerate}
\end{lemma}
\noindent\textit{Proof.}  If $\varphi_k$ is injective, then $\Ker(\varphi_k^{i-k})=\{1\varphi^k\}$, for all $i>k$. By Lemma  \ref{lema nucleos1}, this implies that, for all $i>k$, $\Ker(\varphi^i)=\Ker(\varphi^k)$. So 1 implies 2. Since 2 trivially implies 3, we only have to verify that 3 $\implies$ 1. Suppose that there is some $i>k$ such that  $\Ker(\varphi^k)= \Ker(\varphi^i)$. Then, by Lemma  \ref{lema nucleos1}, it follows that $\varphi_k^{i-k}$ is injective, and so must be $\varphi_k$, since $\Ker(\varphi_k)\leq \Ker(\varphi_k^{i-k})$.
  \qed\\

\subsection{Virtually free groups}

 In \cite{[MS02]}, Myasnikov and Shpilrain prove that there is a periodic orbit of length $k$ in some automorphism of the free group $F_n$ if and only if there is an automorphism of $F_n$ with order $k$. In \cite{[Car22]}, the author shows that the same result holds for endomorphisms. This result allows us to obtain a (computable) bound to the length of every periodic orbit of an endomorphism $\varphi\in \End(F_n)$. 
In \cite[Proposition 5.3]{[Car22]}, the author proves that  such a bound can be obtained for endomorphisms of finitely generated virtually free groups.

  \begin{theorem}
  \label{prop}
 Let $G$ be a f.g. virtually free  group and $\varphi\in \End(G)$. Then, there is a computable constant $C$, depending only on $\varphi$, such that, for all $K=\{h\}\subseteq \Per(\varphi)$, 
 $\spe(K)= [n]$, for some $n\leq C$. Moreover, for such a $K$, $\spe(K)$ is computable.
 \end{theorem}
  \noindent\textit{Proof.}  Let $G$ be a f.g. virtually free  group. From Remark \ref{rem: phik injetivo}, there is a computable constant $k$  such that the restriction of $\varphi$, $\varphi_k:\Img(\varphi^k)\to \Img(\varphi^{k+1})$,  is injective. Also, from \cite[Proposition 5.3]{[Car22]}, there is a computable constant $p$ such that every periodic orbit has length bounded by $p$. Let $C=p+k-1$, $K=\{h\}\subseteq \Per(\varphi)$ and $\pi_h$ denote the period of $h$. Suppose that there is an element $x\in G$ such that $\ord_K(x)=C+1$. Notice that $C+1-\pi_h=p+k-\pi_h\geq k$.
   Then, $$x\varphi^{C+1-\pi_h}\varphi^{\pi_h}=x\varphi^{C+1-\pi_h}\varphi_k^{\pi_h}=h=h\varphi_k^{\pi_h}$$ and so, since $\varphi_k^{\pi_h}$ is injective, then $x\varphi^{C+1-\pi_h}=h$, which contradicts the assumption that $\ord_K(x)=C+1$.
  Hence, we have that $\spe(K)=[n]$, for some $n\leq C$. 
  
  Now, we have that $n\geq \pi_h-1$, since $\ord_K(h\varphi)=\pi_h-1$. We will now prove that, given $\pi_h \leq i\leq C$, we can decide if $i\in \spe(K)$ or not, and that suffices. Fix such $i$. By definition, $x\in G$ is a point of order $i$ if and only if $x\varphi^i=h$ but $x\varphi^j\neq h$ for all $j< i$. This happens exactly when $x\varphi^i=h$ and $x\varphi^{i-\pi_h}\neq h$ because, if $x\varphi^j=h$, for some $j < i$, then, for $r>j$,  $x\varphi^r=h$ if and only if $r=j+k\pi_h$. Write $i=k\pi_h+r$, with $0\leq r<\pi_h$. The set of elements that get mapped to $h$ through $\varphi^i$ (resp. $\varphi^{i-\pi_h}$) are $h\varphi^{\pi_h -r}\Ker(\varphi^{i})$ (resp. $h\varphi^{\pi_h -r}\Ker(\varphi^{i-\pi_h})$).
   It follows that there is an element of order $i$ if and only if $\Ker(\varphi^{i})\setminus\Ker(\varphi^{i-\pi_h})\neq\emptyset$. Since $\Ker(\varphi^{i-\pi_h})\subseteq \Ker(\varphi^{i})$, then $\Ker(\varphi^{i})\setminus\Ker(\varphi^{i-\pi_h})\neq\emptyset$ if and only if $\Ker(\varphi^{i})\neq\Ker(\varphi^{i-\pi_h})$. By Lemma \ref{lem: nucleosphik},
$\Ker(\varphi^{i})=\Ker(\varphi^{i-\pi_h})$ if and only if $\varphi_{i-\pi_h}$ is injective, which we can decide, since it is equivalent to deciding whether $\Img(\varphi^{i-\pi_h})\simeq \Img(\varphi^{i-\pi_h+1})$ or not: if $\varphi_{i-\pi_h}$ is injective, then obviously, $\Img(\varphi^{i-\pi_h})\simeq \Img(\varphi^{i-\pi_h+1})$; we have that  $\Img(\varphi^{i-\pi_h+1})\simeq \faktor{\Img(\varphi^{i-\pi_h})}{\Ker(\varphi_{i-\pi_h})}$ and if  $\Img(\varphi^{i-\pi_h})\simeq \Img(\varphi^{i-\pi_h+1})$, we have that $\Img(\varphi^{i-\pi_h})\simeq \faktor{\Img(\varphi^{i-\pi_h})}{\Ker(\varphi_{i-\pi_h})}$, which implies that $\Ker(\varphi_{i-\pi_h})$ is trivial, because virtually free groups are Hopfian.
 \qed
  \begin{corollary}
  \label{cor: spectrum singleton}
 Let $G$ be a f.g. virtually free  group $K=\{h\}\subseteq G$, and $\varphi\in \End(G)$. Then the following are equivalent:
 \begin{enumerate}
 \item $\spe(K)=\N$
 \item $h\in \varphi^\infty(G)\setminus \Per(\varphi)$
 \end{enumerate}
 \end{corollary}
  \noindent\textit{Proof.} If $\spe(K)=\N$, then, obviously, $h\in \Img(\varphi^k)$, for all $k\in \N$. Also, $h\not\in \Per(\varphi)$ by Proposition \ref{prop}.
  
  Conversely, if $h\in \varphi^\infty(G)\setminus \Per(\varphi)$, then, let $n\in \N$. Since $h\in \Img(\varphi^n)$, there is some $x$ such that $x\varphi^n=h$. Since $h$ is not periodic, we know that $x\varphi^k\neq h$, for every $k< n$, since otherwise we would have $h=x\varphi^n=x\varphi^k\varphi^{n-k}=h\varphi^{n-k}$. Thus $\ord_K(x)=n$ and so $n\in \spe(K)$.
 \qed

   \begin{corollary}
   \label{cor: spectrum periodic}
There exists an algorithm with input a finitely generated virtually free group $G$, an endomorphism $\varphi$ of $G$ and a singleton $K=\{h\}\subseteq G$ and output $\spe(K)$.
 \end{corollary}
 \noindent\textit{Proof.}  Let $G$ be a finitely generated  virtually free group. Then, we can compute $\varphi^\infty(G)$ by Theorem \ref{stable image computable} and test if $h\in \varphi^\infty(G)$. If not, then by successively testing membership of $h$ in $\Img(\varphi^k)$ we compute $m=\max\{k\in \N\mid h\in \Img(\varphi^k)\}$ (might be $0$ if $h\not\in \Img(\varphi)$), and we have that $\spe(K)=[m]_0$. 
 So, assume that $h\in \varphi^\infty(G)$. Since the length of periodic orbits is bounded by a computable constant $p$ and the fixed subgroup of $\varphi$ is computable, then $\Per(\varphi)=\Fix(\varphi^{p!})$ is also computable by \cite[Theorem 4.2]{[Car22]}. If $h\not\in \Per(\varphi)$, then, by Corollary \ref{cor: spectrum singleton}, we have that $\spe(K)=\N$. If, on the other hand, $h\in \Per(\varphi)$, $\spe(K)$ is computable by Theorem \ref{prop}.
\qed\\

It is clear that decidability of $BrP(G)$ implies decidability of $GBrP(G)$ for finite subsets. We will now see that the same occurs for the computability of the spectrum. 
\begin{lemma}
\label{lem: spectrum contained}
Let $G$ be a group, $K\subseteq G$ and $\varphi\in \End(G)$. Then $$\spe(K)\subseteq \bigcup_{y\in K} \spe(\{y\}).$$
\end{lemma}

\noindent\textit{Proof.} Let $n\in \spe(K)$. There is some $x\in G$ such that $x\varphi^n\in K$ and $x\varphi^i\not\in K$ for $i< n$. In particular, $x\varphi^i\neq x\varphi^n$ for $i<n$, so   $n\in \spe(\{x\varphi^n\})\subseteq  \bigcup_{y\in K} \spe(\{y\})$.
\qed

\begin{theorem}\label{main finite}
There exists an algorithm with input a finitely generated virtually free group $G$, an endomorphism $\varphi$ of $G$ and a finite set $K=\{g_1,\ldots, g_k\}\subseteq G$ and output $\spe(K)$.
\end{theorem}
 \noindent\textit{Proof.} 
  We start by proving that $\spe(K)=\N$ if and only if $\spe(\{g_i\})=\N$ for some $i\in\{1,\ldots, k\}$. Suppose that  $\spe(K)=\N$. By Lemma \ref{lem: spectrum contained}, $\spe(K)\subseteq \bigcup_{i=1}^k \spe(\{g_i\})$, and so,  there must be some $i\in\{1,\ldots, k\}$ for which  $\spe(\{g_i\})=\N$.  Now suppose that $\spe(K)=[m]_0$, for some $m\in \N$ and that there is some $i\in \{1,\ldots, k\}$ such that $\spe(\{g_i\})=\N$. Notice that it follows from Theorem \ref{prop} that $g_i$ cannot be a periodic point.
 Put $$I=\{j\in \{1,\ldots, k\}\mid 0<\ord_{\{g_i\}}(g_j)<\infty\}.$$ If $I=\emptyset$, then, we can observe that $\N= \spe(\{g_i\})\subseteq \spe(K)$, since, if $n\in \spe(\{g_i\})$, then there is some $x\in G$ such that $\ord_{\{g_i\}}(x)=n$. Clearly, $\ord_K(x)=n$ since, if there was some $j< n$ for which $x\varphi^j\in K$, then $x\varphi^j=g_r$ for some $r\in\{1,\ldots, k\}$ (and $r\neq i$ since $\ord_{\{g_i\}}(x)=n$) and $$g_r\varphi^{n-j}=x\varphi^j\varphi^{n-j}=x\varphi^n=g_i$$ which, means that $0<\ord_{\{g_i\}}(g_r)\leq n-j<\infty$ and that contradicts the assumption that $I=\emptyset$.
 
 So, suppose that $I\neq\emptyset$ and let $C=\max\limits_{j\in I}\ord_{\{g_i\}}(g_j)$ and $x\in G$ be such that $M=\ord_{\{g_i\}}(x)>m+C$. Since $\ord_{K}(x)\leq m$, there are some $r\in [m]$ and $s\in \{1,\ldots, k\}$ such that $x\varphi^r=g_s$. Since   $$g_s\varphi^{M-r}=x\varphi^r\varphi^{M-r}=x\varphi^M=g_i,$$ then, since $g_i$ is not periodic, $s\in I$ and, by definition of $C$, there must be some $p\leq C$ such that $g_s\varphi^p=g_i$. But, since $M-r>C$, then $M-r-p>0$ and 
  $$g_i\varphi^{M-r-p}=(g_s\varphi^{p})\varphi^{M-r-p}=g_s\varphi^{M-r}=g_i,$$
  which is absurd since $g_i$ is not periodic, by Theorem \ref{prop}.
    
  Therefore, we have proved that $\spe(K)=\N$ if and only if $\spe(\{g_i\})=\N$ for some $i\in\{1,\ldots, k\}$, which is a decidable condition in view of Corollary \ref{cor: spectrum periodic}.
  
  So, assume that, for all $i\in\{1,\ldots, k\}$, $\spe(\{g_i\})=[n_i]_0$, for some $n_i\in \N$. Then, by Lemma \ref{lem: spectrum contained}, $\spe(K)=[m]_0$ for some $m\leq \max_{i\in \{1,\ldots,k\}} n_i$. We will now prove that, given $n\leq  \max_{i\in \{1,\ldots,k\}} n_i$ we can decide whether $n\in \spe(K)$ or not and that suffices. So, fix such an $n.$
  Put
  \begin{align*} 
  \ell=\{i\in \{1,\ldots, k\}\mid \spe(\{g_i\})\subsetneq [n]_0\} \quad \text{ and } \quad L=\{1,\ldots, k\}\setminus \ell.
  \end{align*} 
  These sets are computable by Corollary \ref{cor: spectrum periodic}. If there is some $x\in G$ such that $\ord_K(x)=n$, then $x\varphi^n=g_r$ for some $r\in L$. We will decide, for each $r\in L$, whether there is some $x\in G$ such that 
  \begin{align}
  \label{condition}
  x\varphi^n=g_r \; \wedge \; \forall\, 0\leq i<n ,\;  x\varphi^i\not\in K 
  \end{align}
 or not. Clearly, $n\in \spe(K)$ if and only if there is some $r\in L$ and $x\in G$ for which condition (\ref{condition}) holds. Let  $r\in L$ and put $$I=\{j\in \{1,\ldots, k\}\mid 0<\ord_{\{g_r\}}(g_j)\leq n\}\}.$$ Notice that $I$ is computable since we only have to check if $g_j\varphi^i=g_r$ for $i\leq n_r$.

If $I=\emptyset$, then, it is easy to see that, for any point $x\in G$,  $\ord_{\{g_r\}}(x)=n$  if and only if $x$ satisfies condition (\ref{condition}). This means that there is such an $x$ if and only if $n\in \spe(\{g_r\})$, which is decidable.

So, assume that $I\neq \emptyset$. For each $j\in I$, put $o_j=\ord_{\{g_r\}}(g_j)\leq n$. We now verify if $g_r$ is periodic or not. If $g_r$ is periodic with period $\pi_r$, we add it to $I$ and put $o_r=\pi_r$. If $g_r$ is not periodic, we proceed.

 We claim that, for $x\in G$,  condition (\ref{condition}) holds if and only if
\begin{align}
\label{condition2}
x\varphi^n=g_r\;  \wedge\; \forall j\in I ,\; x\varphi^{n-o_j}\neq g_j.
\end{align}
 Clearly, condition (\ref{condition}) implies  condition (\ref{condition2}). Conversely, if $x\varphi^n=g_r$ but $\ord_K(x)<n$, we will prove that there must be some $j\in I$ such that $x\varphi^{n-o_j}= g_j$. Clearly, there 
 must be some $0\leq s<n$ and some $p\in\{1,\ldots, k\}$ such that $x\varphi^s=g_p$. Take $s$ to be maximal. Then, $$g_r=x\varphi^n=x\varphi^s\varphi^{n-s}=g_p\varphi^{n-s}$$ and so $p\in I$ (notice that this might mean $p=r$ if $g_r$ is periodic). We must have that $n-s\geq o_p$, and so $s\leq n-o_p$.
  If $s=n-o_p$, we are done, since $x\varphi^{n-o_p}=x\varphi^{s}=g_p.$ If $s<n-o_p$, then $x\varphi^{s+o_p}=g_p\varphi^{o_p}=g_r$, which contradicts the maximality of $s$.
  
  Now, we will prove that we can decide if condition (\ref{condition2}) holds for some $r\in L$ and $x\in G$. Compute $a\in g_r\varphi^{-n}$ and, for all $j\in I$, compute some $a_j\in g_j\varphi^{-(n-o_j)}$. We want to decide whether 
  \begin{align}
  \label{condition3}
  a\Ker(\varphi^n)\subseteq \bigcup_{j\in I} a_j\Ker(\varphi^{n-o_j}).
  \end{align}
  Indeed, condition (\ref{condition2}) holds for some $x\in G$ if and only if condition (\ref{condition3}) does not hold. Obviously, (\ref{condition3}) is equivalent to $\Ker(\varphi^n)\subseteq \bigcup_{j\in I} a^{-1}a_j\Ker(\varphi^{n-o_j}).
$
 Now we refine the union, in order to include only the terms for which $a^{-1}a_j\in \Ker(\varphi^n)$ and such that the cosets in the union are all disjoint and we call $I_2$ the new set of indices.
 This can be done, since we can test membership in $\Ker(\varphi^n)$ and, since for $i\geq j$, we have that $\Ker(\varphi^j)\subseteq \Ker(\varphi^i)$, it follows that given two cosets $a\Ker(\varphi^i)$ and $b\Ker(\varphi^j)$, either  $b\Ker(\varphi^j)\subseteq b\Ker(\varphi^i)= a\Ker(\varphi^i)$, in which case we remove the smaller one from the union, or $b\Ker(\varphi^j)\cap a\Ker(\varphi^i)=\emptyset$. Moreover, this can be easily checked. Also, if $a^{-1}a_j\not\in \Ker(\varphi^n)$, then $a^{-1}a_j \Ker(\varphi^{n-o_j})\subseteq a^{-1}a_j \Ker(\varphi^{n})$ is disjoint from $\Ker(\varphi^n)$, and so it can be removed. 

So, we will prove that we can decide if
 \begin{align}
   \label{condition4}
 \Ker(\varphi^n)= \bigcup_{j\in I_2} a^{-1}a_j\Ker(\varphi^{n-o_j}).
\end{align}
 and that concludes the proof.

Assume then that  (\ref{condition4}) holds. Let $I_2=\{j_1,\ldots j_d\}$ where $o_{j_1}\geq \cdots \geq o_{j_d}$.
Suppose that there is some $x\in \Ker(\varphi^{n-o_{j_2}})\setminus \Ker(\varphi^{n-o_{j_1}})$. Then $a^{-1}a_{j_1}x\in \Ker(\varphi^n)$ and so, since $a^{-1}a_{j_1}x\not\in a^{-1}a_{j_1}\Ker(\varphi^{n-o_{j_1}})$, there must be some $1<r\leq d$ such that $a^{-1}a_{j_1}x\in a^{-1}a_{j_r}\Ker(\varphi^{n-o_{j_r}})$. Since $r\geq 2$, then $\Ker(\varphi^{n-o_{j_2}})\subseteq \Ker(\varphi^{n-o_{j_r}})$, thus $x\in \Ker(\varphi^{n-o_{j_r}})$ and so $a^{-1}a_{j_1}\in a^{-1}a_{j_r}\Ker(\varphi^{n-o_{j_r}})$, which contradicts the fact that $a^{-1}a_{j_1}\Ker(\varphi^{n-o_{j_1}})$ and $a^{-1}a_{j_r}\Ker(\varphi^{n-o_{j_r}})$ are disjoint. This implies that $\Ker(\varphi^{n-o_{j_2}})= \Ker(\varphi^{n-o_{j_1}})$ and so, by Lemma \ref{lem: nucleosphik}, 
 $$\Ker(\varphi^{n})=\Ker(\varphi^{n-o_{j_1}})=\cdots = \Ker(\varphi^{n-o_{j_d}}),$$
which, since we are assuming that the union is disjoint implies that $|I_2|=1$. 

So, in order to decide (\ref{condition4}) we only have to verify if $|I_2|=1$ and if $\Ker(\varphi^{n-o_j})=\Ker(\varphi^n)$, for $j\in I_2$,which can be done since it corresponds to checking injectivity of $\varphi_{n-o_j}$, by  Lemma \ref{lem: nucleosphik}.
   \qed\\
\subsection{Virtually abelian groups}

Using similar arguments to the ones in   \cite{[Car22],[MS02]}, we can prove a similar result for free-abelian groups.
\color{black}
\begin{proposition}
If there is an element $x\in \Z^m$ and an endomorphism $Q\in \mc M_{m\times m}(\Z)=\End(\Z^m)$ such that  $x$ is a periodic point of period $k$, then there is a matrix $P\in GL_m(\Z)$ of order $k$. 
\end{proposition}
  \noindent\textit{Proof.}  Let $x\in \Z^m$ and $Q\in \End(\Z^m)=\mc M_{m\times m}(\Z)$ be such that  $xQ^k=x$ and $xQ^i\neq x$, for all $i\in \{1,\ldots, k-1\}$.  We have that $M=Q|_{ Q^\infty(\Z^m)}$ is an automorphism of $ Q^\infty(\Z^m)$ and 
  $x\in \Per(Q)\subseteq Q^\infty(\Z^m)$.
  Let $H=\Fix(M^k)$. Notice that $M|_{H}$ is an automorphism of $H$:
  $H$ is $M$-invariant since, for $h\in \Fix(M^k)$, we have that $hMM^k=hM^kM=hM$, and so $hM\in \Fix (M^k)$; $M|_H$ is surjective, since for $h\in H$, $h=hM^{k-1}M$; and it is injective since we cannot have $hM=1$ and $h=hM^k$ for $h\neq 1.$
  
  By construction, it is clear that  $M|_H$ has order (as an element of $\Aut(H)$) divisible by $k$. Since $x\in H$ has period $k$, then the order of $M|_H$ is exactly $k$. Since $\rank(H)=r\leq m$, then there is an automorphism of $\Z^m$ of order $k$, which can be defined by applying $M|_H$ to the first $r$ generators and the identity to the remaining ones.   
  \qed

\begin{theorem}[\cite{[Min87]}]
The least common multiple of the orders of all finite subgroups of $n\times n$ matrices over $\Q$ is
$$M_n=\prod_p p^{\lfloor \frac{n}{p-1}\rfloor+\lfloor \frac{n}{p(p-1)}\rfloor+\lfloor \frac{n}{p^2(p-1)}\rfloor+\cdots}$$
where $p$ runs over all primes. 
\end{theorem}

\begin{corollary}
\label{cor: bound virt abelian}
There is a computable constant $k$ that bounds the size of the periodic orbits $\Orb_\varphi(x)$, when $\varphi$ runs through $\mc M_{m\times m}(\Z)$ and $x$ runs through $\Z^m$.
\end{corollary}

Mimicking the argument in \cite[Proposition 5.3]{[Car22]}, one can also extend the previous result to f.g. virtually abelian groups.
\begin{corollary}
\label{cor: vabelian periods}
There exists an algorithm with input a finitely generated virtually abelian group $G$ and output a constant $k$ such that, for all $\varphi\in \End(G)$, the infinite ascending chain $$\Fix(\varphi)\subseteq \Fix(\varphi^{2!})\subseteq \Fix(\varphi^{3!})\subseteq \cdots$$
stabilizes after $k$ steps. Equivalently, if $x\in \Evper(\varphi)$, for some endomorphism $\varphi$ of $G$, then the periodic part of the orbit of $x$ has cardinality at most $k$.
\end{corollary}

Using the fact that the cardinality of periodic orbits is bounded and virtually abelian groups are Hopfian, Theorem \ref{prop} and Corollary \ref{cor: spectrum singleton} hold for finitely generated virtually abelian groups, as the same proofs work. However, Theorem \ref{main finite} holds if and only if we can test membership in the stable image of an endomorphism of a finitely generated virtually abelian group. In particular,  
an affirmative answer to  Question \ref{question stable image v abelian} yields a proof of Theorem \ref{main finite} for finitely generated virtually abelian groups.
\color{black}

\section{Recognizable subsets}
We say that $G$ is \emph{index-computable} if, given a finite index subgroup $N$ by its generators, we can compute its index and coset representatives 
 $a_1,\ldots, a_n\in G$ such that $$G= Na_1\cup \cdots \cup Na_n,$$
 where $Na_1=N$.

 Moreover, let $G$ be a group such that the membership problem is decidable for finite index subgroups. Then $G$ must be index-computable, since we can construct the Schreier graph of $G$ with respect to $N$ by succesively testing membership on $N$ until the graph is complete, meaning that, when we multiplying a coset by a generator, the new coset coincides with a previously computed one. In \cite{[Har11]}, Hartung proves that the membership problem for finite index subgroups is decidable for all finitely $L$-presented groups. In particular, all finitely presented groups are index-computable. This also implies that we can compute generators for the intersection of two finite index subgroups, as observed in \cite{[Har11]}.

\begin{lemma}
\label{recnormal} 
Let $G$ be a finitely presented  group and $H$ be a finite index subgroup of $G$. Then, we can compute a fully invariant finite index normal subgroup $N$ of $G$ contained in $H$ and coset representatives $a_1,\ldots, a_n\in G$ such that $G=Na_1\cup \cdots \cup Na_n$ and $b_1,\ldots, b_m\in H$  such that $H=Nb_1\cup \cdots \cup Nb_m$.
\end{lemma}
\noindent\textit{Proof.} We start by computing the index $[G:H]$ and coset representatives $c_1,\ldots, c_r\in G$ such that $G=Hc_1\cup \cdots \cup Hc_r$, with $Hc_1=H$. Let $\rho:G\to S_r$
be the homomorphism induced by the action of $G$ on $\faktor{G}{H}$. We can compute the image of $\rho$ and so write $G$ as a disjoint union $G=\Ker(\rho)d_1\cup\cdots \cup \Ker(\rho)d_s$ for some $d_1,\ldots, d_s\in G$. Since we can test membership in the kernel of $\rho$, we can compute a generating set for $\Ker(\rho)$ by Schreier's Lemma. Obviously $\Ker(\rho)\trianglelefteq_{f.i.} G$ and $\Ker(\rho)\leq H$, since $g\in \Ker(\rho)$ implies that $Hg=Hc_1g=Hc_1=H$.

By Lemma \ref{fullyinvariant}, the intersection $N$ of all  normal subgroups of $G$ of index at most $s$  is a fully invariant finite index normal subgroup, i.e. $[G:N]<\infty$ and $N\varphi\subseteq N$ for all endomorphisms $\varphi\in \End(G)$. Also, since $[G:\Ker(\rho)]=s$ and $\Ker(\rho)\trianglelefteq G$, then $N\leq \Ker(\rho)$.
We will now prove that $N$ is computable. We start by  enumerating all finite groups of cardinality at most $s$. For each such group $K=\{k_1,\ldots, k_t\}$ we enumerate all homomorphisms from $G$ to $K$ by defining images of the generators and checking all the relations. For each homomorphism $\theta:G\to K$, we have that $[G:\Ker(\theta)]=|\Img(\theta)|\leq |K|\leq s$. In fact, all normal subgroups of $G$ of index at most $s$ are of this form. We compute generators for the kernel of each $\theta$, which is possible since we can test membership in $\Ker(\theta)$, which is a finite index subgroup. We can also find $g_1,\ldots,g_{\ell}\in G$ such that 
$$G=\Ker(\theta)g_1\cup\cdots\cup \Ker(\theta)g_{\ell},$$
taking $g_i$ such that $g_i\theta=k_i.$

Hence, we can compute $N$, since it is a finite intersection   of finite index subgroups, and, since $G$ and $H$ are index-computable, we can compute decompositions of $G$ and $H$ as  disjoint unions of the form
  $$G=Na_1\cup \cdots \cup Na_n$$ and $$H=Nb_1\cup \cdots \cup Nb_m.$$
\qed

Now we can prove the main result of this section.

\begin{theorem}
\label{recognizable computability}
Let $G$ be a finitely presented group, $K\in Rec(G)$, $g\in G$ and $\varphi\in \End(G)$. Then $\ord_K(g)$, $\spe(K)$ and the $\varphi$-preorders are computable.
\end{theorem}
\noindent\textit{Proof.} Let $K\in Rec(G)$. Then $K=Ha_1\cup\cdots\cup Ha_r$, for some finite index subgroup $H\leq_{f.i.}G$.
 Then, in view of Lemma \ref{recnormal}, we may  compute some fully invariant finite index normal subgroup $N$ of $G$ contained in $H$ and coset representatives $c_1,\ldots, c_s$ such that $H=Nc_1\cup\cdots \cup Nc_s$. 
Thus, we can compute $b_i\in K$ such that $$K=Nb_1\cup \cdots \cup Nb_k.$$
Since $G$ is index-computable, we can compute a decomposition of $G$ as a disjoint union $G=Nb_1'\cup\cdots\cup Nb_m'$. 
 
 Now,  since $N$ is fully invariant, then $\varphi|_N\in \End(N)$. Testing membership in $N$, compute  $n\in N$ and $\ell\in [m]$ such that $g=nb_\ell'$. 
 Since we can check membership in $N$, we can test membership in $K$. It follows that, if we can decide whether there is some $i\in \N_0$ such that $g\varphi^i\in K$, we can compute  $\ord_K(g)$: if the answer is \texttt{NO}, then $\ord_K(g)=\infty$; if the answer is \texttt{YES}, we compute successive images of $g$ and test membership in $K$ until we obtain a positive answer.
 
 Since $N$ is fully invariant, then $$g\varphi^s=(nb_\ell')\varphi^s=n\varphi^s b_\ell'\varphi^s\in N (b_\ell'\varphi^s).$$
 The mapping $\theta:\faktor{G}{N}\to \faktor{G}{N}$ given by $Na\mapsto N(a\varphi)$ is well-defined: if $Na=Nb$, then $ba^{-1}\in N$ and so $(ba^{-1})\varphi\in N$, which implies that $b\varphi\in N(a\varphi)$, which yields that  $N(b\varphi)=N(a\varphi)$. It is easy to see that it is also a homomorphism. Hence, we can compute the orbit of $Nb_\ell'$ through $\theta$, which is finite.
 
 Now, there is some $i\in \N$ such that $g\varphi^i\in K$ if and only if $Nb_j\in \Orb_\theta(Nb_{\ell}')$ for some $j\in \{1,\ldots, k\}$, which is obviously decidable. Hence, $\ord_K(g)$ is computable. Moreover, this shows that the order is bounded above
by $\left|\faktor{G}{N}\right|=m.$

For a given natural number $n\in \N$, the $\varphi$-preorders of $n$ in $K$ is the subset of elements $x\in G$ such that $x\varphi^n\in K$ but  $x\varphi^k\not\in K$ for all $k\leq n$, which is 
simply $$K\varphi^{-n}\setminus \left(\bigcup_{i=1}^{n-1}  K\varphi^{-i}\right),$$
which is recognizable and computable since it is the intersection of the inverse image of a recognizable subset with the complement of a finite union of inverse images of recognizable subsets.

Given a natural number $n\in \N$, we can decide whether $n\in \spe(K)$ by computing the $\varphi$-preorder of $n$ in $K$. Also, $\spe(K)\subseteq [m]$, and so it is computable.  
\qed

\section{Variations and further work}

Besides Question \ref{question stable image v abelian}, which would allow us to obtain a version of Theorem \ref{main finite} for finitely generated virtually abelian groups, a natural question to ask concerns the quantification of Brinkmann's conjugacy problem.
In \cite{[Bri10]}, Brinkmann also solved the conjugacy version of $BrP(G)$, Brinkmann's Conjugacy Problem (BrCP(G)), for automorphisms of the free group: given $x,y\in G$ and $\varphi\in \Aut(G)$, can we decide whether there is some $n\in \N$ such that $x\varphi^n\sim y$? Decidability of this problem is a crucial step in the proof of decidability of the conjugacy problem for free-by-cyclic groups in \cite{[BMMV06]} and was later generalized to the notion of orbit decidability in \cite{[BMV10]}. The generalized version of this problem is, similar to what happens to $GBrP(G)$, still quite unknown, the more general result being decidability of $GBrCP(G)$ for automorphisms of virtually polycyclic groups in case $K$ is a coset of a finitely generated subgroup (see \cite{[Car22c]}).

 We can analogously define  the \emph{relative $\tilde\varphi$-order of $g$ in $K$}, $\tord_K(g)$, as the smallest nonnegative integer $k$ such that $g\varphi^k$ has a conjugate in $K$  (and if there is no such $k$, we say that   $\tord_K(g)=\infty$).  The \emph{$\tilde\varphi$-spectrum} of a subset $\tspe(K)$ is then the set of relative $\tilde\varphi$-orders of elements in $K$, i.e., $\tspe(K)=\{\tord_K(g)\mid g\in G\}$. Computation of $\tilde\varphi$-spectrum of a subset, seems a complicated, but interesting problem. 

\begin{question}
For a finitely generated virtually free group $G$, an endomorphism $\varphi\in \End(G)$ and a finite subset $K\subseteq G$,  is the $\tilde\varphi$-spectrum  of $K$ computable? 
\end{question}

Another natural question concerns the case where the subset $K$ is a finitely generated subgroup. Solving $GBrP_{f.g.}$, and so computing $\varphi$-orders, is known to be difficult. In  \cite[Corollary 4.5]{[Car22c]}, the author proves that $GBrP_{f.g.}$ is  decidable for automorphisms of finitely generated abelian groups and that is, to the author's knowledge, the most general result on this problem (see \cite[Section 9]{[Ven14]} for additional comments on this topic). Despite being complicated, in the case of automorphisms, the problem admits a seemingly easier formulation.

\subsection{The automorphism case}

 Given an element $g\in G$ and an endomorphism $\varphi\in \End(G)$, we call  \emph{$\varphi$-tail of $g$} to a sequence 
$ (x_{-i})_{i\in \N}$ such that $x_0=g$ and, for all $i>0$,  $x_{-i+1}=x_{-i}\varphi$. We will also refer to finite arrays of this kind as $\varphi$-tail of $g$. So, $(x,x\varphi,\ldots, x\varphi^n)$, where $x\varphi^n=g$ is a $\varphi$-tail of $g$.

Let $t=(x_{-i})_i$ be a $\tail$ of $g$ and $K\subseteq G$. We define 
$$\lambda_K(t):=
\begin{cases}
|t|\quad &\text{if  for all $i>0$, $x_{-i}\not\in K$}\\
\min\{i>0\mid x_{-i}\in K\} \quad &\text{otherwise}\\
\end{cases}
$$
Notice that $\lambda_K(t)\leq |t|$  and that it is possible that $\lambda_K(t)=\infty$.

 In case $\varphi$ is an automorphism, every $\tail$ of  $h$ is contained in the infinite $\tail$ $t=(\ldots, h\varphi^{-2},h\varphi^{-1},h)$.
 \color{black}
 \begin{lemma}
 \label{lem: tails}
  Let $G$ be a group, $K\subseteq G$ and $\varphi\in \End(G)$. Then the following are equivalent:
  \begin{enumerate}
  \item $\spe(K)=[n]_0$
  \item $\max \{\lambda_K(t)\mid g\in K, t \text{ is a $\tail$ of $g$}\}=n+1$
  \end{enumerate}
 \end{lemma}
 \noindent\textit{Proof.} Suppose that $\spe(K)=[n]_0$ and put $S=\{\lambda_K(t)\mid g\in K, t \text{ is a $\tail$ of $g$}\}$. Then, $$K\varphi^{-n-1}\setminus \bigcup_{i=0}^{n} K\varphi^{-i}=\pord_K(n+1)=\emptyset,$$ which means that $K\varphi^{-n-1}\subseteq  \bigcup_{i=0}^{n} K\varphi^{-i}$. 
 So let $g\in K$ and $t=(t_{-i})_i$ be a $\tail$ of $g$.  If $|t|\leq n+1$, then $\lambda_K(t)\leq n+1$. If $|t|>n+1$, by definition of $\varphi$-tail, we have that   $t_{-n-1}\in K\varphi^{-n-1}\subseteq  \bigcup_{i=0}^{n} K\varphi^{-i}$ and so $t_{-n-1}\varphi^k\in K$, for some $k\in [n]_0$, which means that $t_{-n-1+k}\in K$ and so $\lambda_K(t)\leq n+1-k\leq n+1$. 
   Hence, $S$ is bounded above by $n+1$. Clearly, since $\spe(K)=[n]_0$, there is an element $x\in G$ with $\varphi$-order equal to $n$, i.e.,  such that $x\varphi^n\in K$ and $x\varphi^k\not\in K$ for any $0\leq k<n$.
   Then $t=(x,x\varphi,\ldots, x\varphi^n)$ is a $\tail$ of $x\varphi^n\in K$ and $\lambda_K(t)=|t|=n+1.$ So,  $\max \{\lambda_K(t)\mid h\in K, t \text{ is a $\tail$ of $h$}\}=n+1$.
   
   Conversely, suppose that $\max S=n+1$. If there was an element of order $n+1$, arguing as above, we could construct a $\tail$ $t$ of an element $x\in K$ such that   $\lambda_K(t)=n+2$, which proves that $\spe(K)=[k]_0$ for some $k\leq n$. 
   Since $\max S=n+1$, there is a tail $t$ ending in some $g\in K$ such that $\lambda_K(t)=n+1$, and so $t_{-n}$ has order $n$. Therefore $\spe(K)=[n]_0$.
    \qed\\
 
 Naturally,  $\spe(K)=\N$ if and only if $\{\lambda_K(t)\mid h\in K, t \text{ is a $\tail$ of $h$}\}$ is unbounded.

    Since, in the case of automorphisms, every tail ending in an certain element is contained in the unique infinite tail defined by that element, a better version of Lemma \ref{lem: tails} is possible.
    \begin{proposition}
Let $G$ be a finitely generated group, $\varphi\in \Aut(G)$ and $K\subseteq G$. Then, the following are equivalent:
\begin{enumerate}
\item $\spe(K)=[n]_0$;
\item $n$ is the smallest natural number such that 
\begin{align}
\label{condition for autos spectrum}
K\subseteq \bigcup_{k=1}^{n+1} K\varphi^k .
\end{align}

\end{enumerate}
\end{proposition}
\noindent\textit{Proof.} Suppose that $\spe(K)=[n]_0$ and let $g\in K$. Then, consider the $\tail$ of $g$ of length $n+2$: $t=(g\varphi^{-n-1},g\varphi^{-n},\cdots, g\varphi^{-1},g)$. By Lemma \ref{lem: tails}, $\lambda_K(t)\leq n+1$, which means that
there is some $1\leq k\leq n+1$ such that $g\varphi^{-k}\in K$ and so $g\in K\varphi^k\subseteq  \bigcup_{k=1}^{n+1} K\varphi^k$. 
Now, suppose that $K\subseteq \bigcup_{k=1}^{i} K\varphi^k $ for some $i\leq n$ and take $x\in G$ such that $\ord_K(x)=n$. Then, $x\varphi^{n}$ belongs to $K\varphi^k$ for some $1\leq k \leq i\leq n$ which means that $x\varphi^{n-k}\in K$, which contradicts the fact that $\ord_K(x)=n$.

Conversely, suppose that $n$ is the smallest natural number such that  $K\subseteq \bigcup_{k=1}^{n+1} K\varphi^k $. Let $x\in K\setminus \bigcup_{k=1}^{n} K\varphi^k$. Then, $x\in K\varphi^{n+1}\setminus \bigcup_{k=1}^{n} K\varphi^k$. Hence, the $\tail$ of $x$ with length $n+2$ starts and ends in $K$ and has no other element in $K$, which means that $\ord_K(x\varphi^{-n})=n$. Also,
$$\pord_K(n+1)=K\varphi^{-n-1}\setminus \bigcup_{i=0}^{n} K\varphi^{-i}=\emptyset,$$
because  $K\subseteq \bigcup_{k=1}^{n+1} K\varphi^k $ and so, by \ref{defspecpord}, we have that $n+1\not\in \spe(K)$.
\qed\\

Condition (\ref{condition for autos spectrum}) looks more tractable for automorphisms of well-behaved classes of groups and might lead to answering affirmatively the following question with respect to some class $\mc C$.
\begin{question}
Let $G$ be a finitely generated group in $\mc C$, $\varphi\in \Aut(G)$ and $K\leq_{f.g} G$. Can we compute $\spe(K)$?
\end{question}
\color{black}

\section*{Acknowledgments}
The author is grateful to Pedro Silva for fruitful discussions of these topics, which improved
the paper. This work was partially supported by the grant SFRH/BD/145313/2019 funded by the
 FCT - Fundação para a Ciência e a Tecnologia, I.P., and by national funds through the FCT - Fundação para a Ciência e a Tecnologia, I.P., under the scope of the projects UIDB/00297/2020 and UIDP/00297/2020 (Center for Mathematics and Applications).

\bibliographystyle{plain}
\bibliography{Bibliografia}

\begin{thebibliography}{10}

\bibitem{[Ber79]}
J.~Berstel.
\newblock {\em Transductions and Context-free Languages}.
\newblock Teubner, Stuttgart, 1979.

\bibitem{[BMMV06]}
O.~Bogopolski, A.~Martino, O.~Maslakova, and E.~Ventura.
\newblock The conjugacy problem is solvable in free-by-cyclic groups.
\newblock {\em Bull. London Math. Soc.}, 38:787--794, 2006.

\bibitem{[BMV10]}
O.~Bogopolski, A.~Martino, and E.~Ventura.
\newblock Orbit decidability and the conjugacy problem for some extensions of
  groups.
\newblock {\em Transactions of the American Mathematical Society},
  362(4):2003--2036, 2010.

\bibitem{[Bri10]}
P.~Brinkmann.
\newblock Detecting automorphic orbits in free groups.
\newblock {\em Journal of Algebra}, 324(5):1083--1097, 2010.

\bibitem{[Car22]}
A.~Carvalho.
\newblock Eventually fixed points of endomorphisms of virtually free groups.
\newblock {\em preprint, arXiv: 2204.04543}, 2022.

\bibitem{[Car22c]}
A.~Carvalho.
\newblock On generalized conjugacy and some related problems.
\newblock {\em Comm. Algebra}, 51(8):3528--3542, 2023.

\bibitem{[CD23]}
A.~Carvalho and J.~Delgado.
\newblock Homomorphisms of free-abelian-by-free groups.
\newblock {\em in preparation}.

\bibitem{[DVZ22]}
J.~Delgado, E.~Ventura, and A.~Zakharov.
\newblock Relative order and spectrum in free and related groups.
\newblock {\em Commun. Contemp. Math.}, to appear.

\bibitem{[Har11]}
R.~Hartung.
\newblock Coset enumeration for certain infinitely presented groups.
\newblock {\em Internat. J. Algebra Comput.}, 21(8):1369--1380, 2011.

\bibitem{[IT89]}
W.~Imrich and E.~Turner.
\newblock Endomorphisms of free groups and their fixed points.
\newblock {\em Math. Proc. Cambridge Philos. Soc.}, 105(3):421--422, 1989.

\bibitem{[KL86]}
R.~Kannan and R.~Lipton.
\newblock Polynomial-time algorithm for the orbit problem.
\newblock {\em J. Assoc. Comput. Mach.}, 33(4), 1986.

\bibitem{[Log22]}
A.~D. Logan.
\newblock The conjugacy problem for ascending {HNN}-extensions of free groups.
\newblock {\em preprint, ar{X}iv: 2209.04357}, 2022.

\bibitem{[Min87]}
H.~Minkowski.
\newblock Zur theorie der positiven quadratischen formen.
\newblock {\em J. Reine Angew. Math.}, 101:196--202, 1887.

\bibitem{[Mut21]}
J.~P. Mutanguha.
\newblock Constructing stable images.
\newblock preprint, available at https://mutanguha.com/pdfs/relimmalgo.pdf.

\bibitem{[MS02]}
A.~G. Myasnikov and V.~Shpilrain.
\newblock Automorphic orbits in free groups.
\newblock {\em Journal of Algebra}, 269:18--27, 2002.

\bibitem{[Ven14]}
E.~Ventura.
\newblock Group-theoretic orbit decidability.
\newblock {\em Groups Complexity Cryptology}, 6(2):133--148, 2014.

\end{thebibliography}

 \end{document}